\documentclass[11pt]{article}
\usepackage[margin=1in]{geometry}
\usepackage{amsmath,amssymb,amsthm,url}
\usepackage{needspace}
\usepackage[hidelinks]{hyperref}
\newtheorem{theorem}{Theorem}
\newtheorem{lemma}[theorem]{Lemma}
\newtheorem{corollary}[theorem]{Corollary}
\theoremstyle{definition}
\newtheorem{example}[theorem]{Example}
\newtheorem{question}{Question}
\newtheorem{definition}{Definition}
\theoremstyle{remark}

\newcommand{\Ftwo}{\mathbb F_2}
\title{Real Principal-Minor Detectors: Normal Forms and Regular Parity Lifts}
\author{Qingying Deng \qquad Chaoyang Zhang\thanks{Corresponding author.
Email: \href{mailto:202521511267@smail.xtu.edu.cn}{\nolinkurl{202521511267@smail.xtu.edu.cn}}.}\\[0.5em]
\small School of Mathematics and Computational Science, Xiangtan University,
P.R. China}
\date{}
\begin{document}
\maketitle
\begin{abstract}
Let \(A\) be a symmetric matrix over \(\Ftwo\) on a finite index set,
and let \(t\) be the column of diagonal entries of \(A\), viewed over
\(\mathbb R\). A real detector of \(A\) is a real matrix whose principal
minors equal the integer \(0/1\) representatives of the corresponding
principal minors of \(A\). We prove that \(A\) has a real detector if
and only if it has one of the form \(K+tt^T\), with \(K\) skew-symmetric
and \(\{0,\pm1\}\)-valued. These conditions are equivalent to the
existence of a principally unimodular skew-symmetric signing of the
bordered lift matrix of \(A\), and to regularity of the parity lift of
its delta-matroid. If the graph on the nonempty index set of \(A\),
with edges given by its nonzero off-diagonal entries, is connected,
every real detector of \(A\) is diagonally similar to a detector of
this form. The proof uses cut-transposes to preserve the
prescribed symmetric part up to diagonal similarity, then applies
feasible singleton pivots to a minimal counterexample. For ribbon
bouquets, the characterization proves Whistler's conjecture. For every
\(A\) with a real detector, the feasible-set generating polynomial is
Hurwitz stable, and the polynomial counting feasible sets of sizes
\(2i-1\) or \(2i\) by \(i\) has only strictly negative real zeros and
ultra-log-concave coefficients.
\end{abstract}

\section{Introduction}\label{sec:introduction}

A spanning quasi-tree of a connected ribbon graph is a spanning
ribbon subgraph with one boundary component. For a plane ribbon graph,
its spanning quasi-trees are precisely its spanning trees. The term
was introduced by Dasbach, Futer, Kalfagianni, Lin, and Stoltzfus
\cite{DFKLS}. Matrix--Quasi-tree theorems give determinantal formulas
for quasi-tree enumeration in classes of ribbon graphs \cite{MMN,DK}.

For a ribbon bouquet, the edge sets of spanning quasi-trees are exactly
the subsets indexing nonsingular principal submatrices of its framed
intersection matrix over \(\Ftwo\) \cite[Proposition~2.21]{DK}.
These edge sets are also the feasible sets of its ribbon-graphic
delta-matroid; delta-matroids were introduced by Bouchet \cite{Bouchet87}.
The matrix description leads to a principal-minor question: when can
the resulting \(0/1\) data be realized as the principal minors of a
real matrix?

For orientable ribbon graphs, Merino, Moffatt, and Noble proved a weighted
Matrix--Quasi-tree Theorem and Hurwitz stability
\cite[Theorems~6.8--6.10]{MMN}. Ding and Kim introduced
pseudo-orientable ribbon graphs and extended these conclusions to that
class through regular delta-matroids
\cite[Theorem~1.2 and Section~4]{DK}. Their work also isolates the
underlying principal-minor problem: a real matrix detects a bouquet
when its principal minors are exactly the quasi-tree indicators
\cite[Section~5.1]{DK}. Expanding the determinant of the identity plus
a detector counts the spanning quasi-trees.

Whistler used the framed intersection matrix of a bouquet to construct
a detector from a principally unimodular (PU) skew-symmetric signing
of its bordered lift matrix, defined below. He conjectured that every detectable bouquet admits such a signing
\cite[Observation~6.1 and Conjecture~6.2]{Whistler}. Whistler also observed that,
for a bouquet, PU-orientability of the bordered lift matrix is equivalent to
regularity of the parity lift \cite[Observation~6.3]{Whistler}.

We prove the implication from real detectability to PU-orientability
for every finite symmetric binary matrix. The proof also yields a
detector with a prescribed symmetric part.

All index sets are finite. For a matrix \(M\) indexed by \(E\) and
\(X\subseteq E\), let \(M[X]\) denote the principal submatrix on \(X\).
We set \(\det M[\varnothing]=1\).
Write \(M[U,V]\) for the submatrix on row set \(U\) and column set
\(V\), and \(v[X]\) for the restriction of a column vector \(v\) to \(X\).
Write \(\triangle\) for symmetric difference. For a symmetric binary
matrix \(A\in\Ftwo^{E\times E}\), let \(t=\operatorname{diag}A\) be the
column of its diagonal entries. We write \(\det\nolimits_2\) for the
determinant over \(\Ftwo\).
In real matrix identities, counting sums, and generating polynomials,
binary scalars and vectors are interpreted using their integer \(0/1\)
representatives.
The off-diagonal support graph of \(A\) has vertex set \(E\) and
edge \(ij\) when \(i\ne j\) and \(A_{ij}=1\).
A symmetric binary matrix with zero diagonal is called alternating.

\begin{definition}\label{def:bordered-lift}
Let \(p\notin E\). Following \cite[Section~6]{Whistler}, define the
bordered lift matrix of \(A\) by
\begin{equation}\label{eq:general-lift}
 A^\sharp=
 \begin{pmatrix}A+tt^T&t\\t^T&0\end{pmatrix}
 \in\Ftwo^{(E\cup\{p\})\times(E\cup\{p\})}.
\end{equation}
The matrix \(A^\sharp\) is alternating because the diagonal of
\(A+tt^T\) is zero over \(\Ftwo\).
\end{definition}

\begin{definition}\label{def:PU}
A matrix is principally unimodular (PU) if its principal minors lie in \(\{0,1,-1\}\). For an alternating binary matrix \(H\)
on an index set \(V\), a PU signing of \(H\) is a PU skew-symmetric
matrix \(L\in\{0,1,-1\}^{V\times V}\) with \(L\bmod2=H\).
We call a PU signing of \(A^\sharp\) a PU-orientation of \(A^\sharp\),
following \cite[Section~6]{Whistler}.
We write \(\widetilde{A^\sharp}\) for a chosen PU-orientation of
\(A^\sharp\).
\end{definition}

\begin{definition}\label{def:detector}
A real detector of \(A\), extending the bouquet terminology of
\cite[Section~5.1]{DK}, is a matrix \(M\in\mathbb R^{E\times E}\)
such that
\begin{equation}\label{eq:detector}
 \det M[X]=\det\nolimits_2 A[X]\qquad(X\subseteq E).
\end{equation}
A special detector is a detector of the form \(K+tt^T\), where
\(K\) is integer skew-symmetric with entries in \(\{0,1,-1\}\).
\end{definition}

A real skew-symmetric matrix has determinant zero in odd order;
in even order its determinant is the square of its Pfaffian
\cite[Section~I.A]{Wimmer}.
Thus a skew-symmetric matrix is PU precisely when all its
principal minors belong to \(\{0,1\}\).

\begin{theorem}\label{thm:binary-detector}
\label{thm:general-bouquet}
Let \(E\) be finite, let \(A=A^T\in\Ftwo^{E\times E}\), and put
\(t=\operatorname{diag}A\). Choose \(p\notin E\), and let
\(A^\sharp\) be the bordered lift matrix in
Definition~\ref{def:bordered-lift}.
The following conditions are equivalent.
\begin{enumerate}
\item There exists \(M\in\mathbb R^{E\times E}\) such that
 \(\det M[X]=\det\nolimits_2 A[X]\) for every \(X\subseteq E\).
\item There exists an integer skew-symmetric
 \(K\in\{0,1,-1\}^{E\times E}\) such that
 \(\det(K+tt^T)[X]=\det\nolimits_2 A[X]\) for every \(X\subseteq E\).
\item The matrix \(A^\sharp\) admits a PU-orientation.
\end{enumerate}
\end{theorem}

Corollary~\ref{cor:binary-regularity} gives the equivalent regularity
criterion. Corollary~\ref{cor:general-connected-normal} shows that on
nonempty connected support, every real detector is diagonally similar
to a special detector.

Corollary~\ref{cor:bouquet-main} gives the formulation for bouquets.
The topological definitions and the principal-minor description of
quasi-trees are given in Section~\ref{sec:bouquet-applications}.
Chatterjee, Ghosh, Gurjar, and Raj characterized
principal-minor equivalence of irreducible matrices by cut-transposes
and diagonal equivalence \cite[Theorem~1.1]{CGGR}.
Lemma~\ref{lem:general-fibre} uses that characterization to preserve
a prescribed symmetric part up to diagonal similarity.
Lemma~\ref{lem:zero-diagonal} supplies the zero-diagonal case directly
from the detector equations. The proof in
Section~\ref{sec:general-binary} combines these inputs with proper
restrictions and feasible singleton pivots.

Section~\ref{sec:generating-polynomials} gives polynomial consequences
for a finite symmetric binary matrix that has a real detector. Its
feasible-set generating polynomial is Hurwitz stable. The polynomial
whose coefficient of \(z^i\) counts feasible sets of size \(2i-1\) or
\(2i\) has only negative real zeros.

Ding and Kim proved Hurwitz stability for the principal-minor
polynomial of \(K+tt^T\), where \(K\) is real skew-symmetric and
\(t\) is a real vector \cite[Lemma~4.4]{DK}. Our stability consequence
starts from the existence of an arbitrary real detector and uses
Theorem~\ref{thm:binary-detector} to obtain a detector of the required form.


\section{Preliminaries}\label{sec:preliminaries}

\subsection{Principal minors of orders one and two}

Let \(A=A^T\in\Ftwo^{E\times E}\) and let
\(t=\operatorname{diag}A\).
Singleton and two-element equations for a real detector \(M\) give
\begin{equation}\label{eq:detector-pairs}
 M_{ii}=t_i,\qquad
 M_{ij}M_{ji}=
 \begin{cases}
  2t_it_j-1,&A_{ij}=1,\\
  0,&A_{ij}=0
 \end{cases}
 \quad(i\ne j).
\end{equation}
For \(i\ne j\), the detector property gives
\[
 \det M[\{i,j\}]=t_it_j-M_{ij}M_{ji}.
\]
Over \(\Ftwo\),
\[
 \det\nolimits_2 A[\{i,j\}]=t_it_j+A_{ij}.
\]
Its integer representative is \(t_it_j\) when \(A_{ij}=0\), and
\(1-t_it_j\) when \(A_{ij}=1\). Comparing these values with
\(\det M[\{i,j\}]\) yields \eqref{eq:detector-pairs}.
For any real skew-symmetric \(K\) indexed by \(E\), set \(R=K+tt^T\).
For distinct \(i,j\in E\), the identities \(t_i,t_j\in\{0,1\}\) give
\begin{equation}\label{eq:rankone-pair}
\begin{aligned}
 \det R[\{i,j\}]
 &=t_it_j-(t_it_j+K_{ij})(t_it_j-K_{ij})\\
 &=t_it_j-(t_it_j)^2+K_{ij}^2\\
 &=K_{ij}^2.
\end{aligned}
\end{equation}
Now suppose \(R\) is a special detector of \(A\). Then \(K_{ij}^2\) is the integer representative of
\(\det\nolimits_2 A[\{i,j\}]\). Hence
\[
 K_{ij}^2\equiv t_it_j+A_{ij}\pmod2.
\]
Since \(K_{ij}\in\{0,\pm1\}\), we have \(K_{ij}\equiv K_{ij}^2\pmod2\).
Thus, for distinct \(i,j\), \(K_{ij}\equiv A_{ij}+t_it_j\pmod2\).
The diagonal of \(K\) is zero, and the diagonal of \(A+tt^T\)
is also zero over \(\Ftwo\), since \(t_i=A_{ii}\) and \(t_i^2=t_i\).
Therefore,
\begin{equation}\label{eq:special-reduction}
 K\bmod2=A+tt^T.
\end{equation}

\subsection{Principal pivot transforms}

Let \(H\) be a square matrix over a field indexed by \(E\), and let
\(F\subseteq E\) satisfy that \(H[F]\) is invertible. In the block order
\(F,E\setminus F\), its principal pivot transform \(H*F\) is
\[
 \begin{pmatrix}P&Q\\R&S\end{pmatrix}*F
 =\begin{pmatrix}P^{-1}&-P^{-1}Q\\
                  RP^{-1}&S-RP^{-1}Q\end{pmatrix}.
\]
The transform at the empty set is the identity.
Multiplying the block formulas gives \((H*F)*F=H\).
Tucker's principal-minor identity, in this convention, is
\begin{equation}\label{eq:ppt}
 \det((H*F)[X])=\frac{\det H[F\triangle X]}{\det H[F]}
 \qquad(X\subseteq E);
\end{equation}
see \cite[Theorem~8]{CG}, where the identity is attributed to Tucker.
Its application to detectable bouquets also appears in
\cite[proof of Lemma~5.2]{DK}.
For symmetric binary \(A\), a set \(F\) is feasible when
\(\det\nolimits_2 A[F]=1\).
In particular, the singleton pivot \(A*\{i\}\) is feasible exactly when \(A_{ii}=1\).
The binary block formula is symmetric, since minus equals plus
in \(\Ftwo\). Thus \(A*F\) is again binary symmetric.
For a real skew-symmetric matrix, write \(R=-Q^T\), \(P^T=-P\),
and \(S^T=-S\) in the displayed block formula.
Then \((P^{-1})^T=-P^{-1}\) and
\[
 (S+Q^TP^{-1}Q)^T=-S-Q^TP^{-1}Q.
\]
The transpose of the upper right block \(-P^{-1}Q\) is
\(Q^TP^{-1}\), the negative of the lower left block.
Thus the principal pivot transform of a real skew-symmetric matrix is
skew-symmetric as well.

\subsection{Delta-matroids and parity lifts}

\begin{definition}
A delta-matroid on \(E\) is a nonempty family of subsets, called
feasible sets, satisfying the symmetric exchange axiom \cite{Bouchet87,BCG}: for feasible \(X,Y\)
and \(e\in X\triangle Y\), there is \(f\in X\triangle Y\) such that
\(X\triangle\{e,f\}\) is feasible. The choice \(f=e\) is allowed.
It is even if all feasible sets have the same parity.
Twisting by \(U\subseteq E\) replaces each feasible set \(X\) by
\(X\triangle U\).
\end{definition}

Write \(\mathcal F(D)\) for the feasible-set family of a delta-matroid \(D\).

For symmetric binary \(A\), its nonsingular principal subsets form
the delta-matroid \(D(A)\), by \cite[Theorem~2.15]{DK}.
For a real skew-symmetric matrix \(L\), the same theorem gives a
delta-matroid \(D(L)\). It is even because every odd-order principal
minor of \(L\) vanishes.
Equation~\eqref{eq:ppt} gives \(D(H*F)=D(H)*F\) for a binary
symmetric or real skew-symmetric matrix \(H\) whenever the pivot
is defined.

For the fixed \(p\notin E\), define \(\alpha:2^E\to 2^{E\cup\{p\}}\) by
\[
 \alpha(X)=
 \begin{cases}X,&|X|\text{ even},\\
 X\cup\{p\},&|X|\text{ odd}.
 \end{cases}
\]
The parity lift \(\widehat{D(A)}\) has feasible sets
\(\{\alpha(X):\det\nolimits_2 A[X]=1\}\), as in \cite[Definition~2.26]{DK}.
It is represented over \(\Ftwo\) by the bordered lift matrix
\(A^\sharp\) of Definition~\ref{def:bordered-lift}.
In particular,
\begin{equation}\label{eq:lift-minors}
 \det\nolimits_2 A^\sharp[\alpha(X)]=\det\nolimits_2 A[X]\qquad(X\subseteq E);
\end{equation}
this is the matrix form of the lift in
\cite[Lemma~2.29 and Proposition~2.32]{DK}.
Every even-sized subset of \(E\cup\{p\}\) is \(\alpha(X)\) for
exactly one \(X\subseteq E\). All odd principal minors of the
alternating binary matrix \(A^\sharp\) vanish.

\begin{definition}
An even delta-matroid is regular if it has a representation
\(D(L)*F\), where \(L\) is a real PU skew-symmetric matrix and
\(F\) is a subset of its index set.
\end{definition}

This is the PU representation definition of regularity;
its field-theoretic equivalents are given in \cite[Theorem~1.3]{BCG}.

\begin{lemma}\label{lem:regular-signing}
Let \(H\) be an alternating binary matrix on a finite set \(V\).
Then \(D(H)\) is regular if and only if \(H\) has a PU signing.
\end{lemma}

\begin{proof}
Let \(L\) be a PU signing of \(H\).
Its principal minors lie in \(\{0,1\}\) by principal unimodularity
and skew-symmetry.
Since \(L\bmod2=H\), the determinant polynomial gives
\[
 \det L[Y]\equiv\det\nolimits_2 H[Y]\pmod2
 \qquad(Y\subseteq V).
\]
The left side lies in \(\{0,1\}\) and the right side is its integer
representative, so the two integers are equal.
Thus \(D(L)=D(H)\).
Twisting by \(\varnothing\) is the identity, and this is the
representation \(D(L)*\varnothing\). Hence \(D(H)\) is regular.

Conversely, \(\det H[\varnothing]=1\), so \(\varnothing\) is feasible
in \(D(H)\).
Write \(D(H)=D(L)*F\), where \(L\) is a real PU skew-symmetric matrix and \(F\subseteq V\).
A set \(Z\) is feasible in \(D(L)*F\) precisely when \(Z\triangle F\) is
feasible in \(D(L)\). Taking \(Z=\varnothing\) shows that \(F\) is
feasible in \(D(L)\).
Principal minors of \(L\) lie in \(\{0,1\}\), so
\(\det L[F]=1\).
Set \(L'=L*F\).
The skew-symmetry calculation above shows that \(L'\) is
skew-symmetric.
Equation~\eqref{eq:ppt} and \(\det L[F]=1\) give
\[
 \det L'[X]=\det L[F\triangle X]\in\{0,1\}
 \qquad(X\subseteq V),
\]
and therefore \(D(L')=D(L)*F=D(H)\).
The diagonal of \(L'\) is zero. For distinct \(i,j\),
\[
 \det L'[\{i,j\}]=(L'_{ij})^2\in\{0,1\},
\]
so \(L'_{ij}\in\{0,1,-1\}\).
Over \(\Ftwo\), the same order-two determinant is
\(\det\nolimits_2 H[\{i,j\}]=H_{ij}\).
Thus \(\{i,j\}\) is feasible in \(D(H)\) precisely when \(H_{ij}=1\),
and feasible in \(D(L')\) precisely when \(L'_{ij}\ne0\).
The two families of pairs agree. Since \(-1\equiv1\pmod2\) and both
diagonals are zero, \(L'\bmod2=H\).
Therefore \(L'\) is a PU signing of \(H\).
\end{proof}

\subsection{Principal-minor equivalence}

Two real matrices on the same index set are principal-minor
equivalent if all corresponding principal minors agree. A
principal-minor fibre is an equivalence class under this relation.
They are diagonally similar if one is \(D^{-1}MD\) for an invertible
real diagonal matrix \(D\). Such a similarity preserves principal
minors because
\((D^{-1}MD)[X]=D[X]^{-1}M[X]D[X]\).
A real matrix is irreducible when the directed graph with arcs
\(i\to j\) for its nonzero off-diagonal entries is strongly
connected. The one-vertex graph is included.
These principal-minor equivalence and irreducibility conventions
are those of \cite[Sections~2.1--2.2]{CGGR}.

\section{Borders, feasible pivots, and direct sums}\label{sec:operations}\label{sec:pivots}

Throughout this section, \(A=A^T\in\Ftwo^{E\times E}\) has a nonempty
index set \(E\). We use \(t\), \(p\), and \(\alpha\) as defined above.

\subsection{The border identity and PU-orientations}

Let \(K\) be any real skew-symmetric matrix on \(E\), let \(u\) be a
real column indexed by \(E\), and put
\(L=\left(\begin{smallmatrix}K&u\\-u^T&0\end{smallmatrix}\right)\), with
the last index \(p\). The border identity is
\begin{equation}\label{eq:border-parity}
 \det\bigl((K+uu^T)[X]\bigr)=\det L[\alpha(X)]\qquad(X\subseteq E).
\end{equation}
This is \cite[Lemma~2.29]{DK} over \(\mathbb R\). Whistler uses it
to construct a detector from a PU-orientation
\cite[Observation~6.1]{Whistler}.

\begin{lemma}\label{lem:general-border}
Let \(A,t\) be as above. Suppose that \(K\) is an integer
skew-symmetric \(E\times E\) matrix with entries in \(\{0,1,-1\}\),
that \(K\bmod2=A+tt^T\), and that \(R=K+tt^T\) detects \(A\).
Then
\[
 \widetilde{A^\sharp}:=
 \begin{pmatrix}K&t\\-t^T&0\end{pmatrix}
\]
is a PU-orientation of \(A^\sharp\). More precisely,
\[
 \det\bigl(\widetilde{A^\sharp}[\alpha(X)]\bigr)
 =\det\nolimits_2 A[X]\qquad(X\subseteq E),
\]
and every odd-order principal minor of \(\widetilde{A^\sharp}\) is zero.
\end{lemma}

\begin{proof}
The matrix \(\widetilde{A^\sharp}\) is
skew-symmetric, since \(K^T=-K\) and the off-diagonal blocks are
\(t\) and \(-t^T\). The entries of \(K\) lie in \(\{0,1,-1\}\)
and those of \(t\) are \(0\) or \(1\), so every entry of
\(\widetilde{A^\sharp}\) lies in \(\{0,1,-1\}\).

Equation~\eqref{eq:border-parity} and the detector property give
\[
 \det\bigl(\widetilde{A^\sharp}[\alpha(X)]\bigr)
 =\det R[X]=\det\nolimits_2 A[X]\qquad(X\subseteq E).
\]
The map \(\alpha\) covers all even subsets, and the odd-order principal
minors vanish by skew-symmetry.
Thus every principal minor of \(\widetilde{A^\sharp}\) lies in \(\{0,1\}\).
By the hypothesis \(K\bmod2=A+tt^T\) and the congruence
\(-t^T\equiv t^T\pmod2\),
\[
 \widetilde{A^\sharp}\bmod2
 =\begin{pmatrix}A+tt^T&t\\t^T&0\end{pmatrix}
 =A^\sharp.
\]
Thus \(\widetilde{A^\sharp}\) is a PU-orientation of \(A^\sharp\).
\end{proof}

\subsection{Invariance under feasible pivots}

\begin{lemma}\label{lem:pivot-invariance}
Let \(F\subseteq E\) satisfy \(\det\nolimits_2 A[F]=1\), and use the
same border index \(p\) for \(A^\sharp\) and \((A*F)^\sharp\).
Then \(A\) has a real detector if and only if \(A*F\) has a real
detector. Moreover, \(A^\sharp\) has a PU signing if
and only if \((A*F)^\sharp\) has a PU signing.
\end{lemma}

\begin{proof}
The binary transform \(A*F\) is symmetric.
Suppose that \(M\) detects \(A\).
Then \(\det M[F]=\det\nolimits_2 A[F]=1\).
For every \(X\subseteq E\), Tucker's identity~\eqref{eq:ppt} over
\(\mathbb R\), the detection identity for \(M\), and
Tucker's identity over \(\Ftwo\) give
\[
 \det\bigl((M*F)[X]\bigr)
 =\det M[F\triangle X]
 =\det\nolimits_2 A[F\triangle X]
 =\det\nolimits_2\bigl((A*F)[X]\bigr).
\]
Thus \(M*F\) detects \(A*F\).

The set \(F\) is also feasible for \(A*F\), because
\[
 \det\nolimits_2\bigl((A*F)[F]\bigr)
 =
 \frac{\det\nolimits_2 A[\varnothing]}
      {\det\nolimits_2 A[F]}
 =1.
\]
If \(N\) detects \(A*F\), the same argument applied to \(A*F\) and
\(F\) shows that \(N*F\) detects
\[
 (A*F)*F=A.
\]
This proves the first equivalence.

We next compare the parity lifts. For every \(X\subseteq E\), the
sets
\[
 \alpha(X\triangle F)
 \quad\text{and}\quad
 \alpha(X)\triangle\alpha(F)
\]
have the same intersection with \(E\), namely \(X\triangle F\).
Moreover, \(p\) belongs to either set exactly when
\(|X\triangle F|\) is odd. Indeed,
\[
 |X\triangle F|
 \equiv |X|+|F|\pmod2,
\]
so this occurs exactly when \(|X|\) and \(|F|\) have different
parities. Hence
\[
 \alpha(X\triangle F)
 =
 \alpha(X)\triangle\alpha(F).
\]

By the parity-lift identity~\eqref{eq:lift-minors} and
Tucker's identity~\eqref{eq:ppt},
\[
\begin{aligned}
 \mathcal F\bigl(D((A*F)^\sharp)\bigr)
 &=
 \{\alpha(Y):
   \det\nolimits_2(A*F)[Y]=1\}\\
 &=
 \{\alpha(Y):
   \det\nolimits_2 A[F\triangle Y]=1\}\\
 &=
 \{\alpha(X\triangle F):
   X\in\mathcal F(D(A))\}\\
 &=
 \{\alpha(X)\triangle\alpha(F):
   X\in\mathcal F(D(A))\}\\
 &=
 \mathcal F\bigl(D(A^\sharp)*\alpha(F)\bigr).
\end{aligned}
\]
Therefore
\begin{equation}\label{eq:lift-twist}
 D((A*F)^\sharp)=D(A^\sharp)*\alpha(F).
\end{equation}

Regularity is invariant under twisting. More precisely, if
\(D=D(L)*U\) for a real PU skew-symmetric matrix \(L\), then
\[
 D*\alpha(F)=D(L)*(U\triangle\alpha(F)),
\]
so \(D*\alpha(F)\) is regular. Conversely,
\((D*\alpha(F))*\alpha(F)=D\), and the same argument shows that
regularity of \(D*\alpha(F)\) implies regularity of \(D\).
By \eqref{eq:lift-twist}, regularity of \(D(A^\sharp)\) is equivalent
to regularity of \(D((A*F)^\sharp)\).
The required equivalence of PU signings now follows
from Lemma~\ref{lem:regular-signing}.
\end{proof}

\subsection{Direct sums of special detectors}

Let \(s\ge1\), let \(E_1,\ldots,E_s\) be pairwise disjoint finite
nonempty sets, and put \(E=E_1\cup\cdots\cup E_s\).
For each \(j\), let \(A_j=A_j^T\in\Ftwo^{E_j\times E_j}\).
Write
\[
 t_j=\operatorname{diag}A_j,\qquad
 A=\operatorname{diag}(A_1,\ldots,A_s),\qquad
 t=(t_1^T,\ldots,t_s^T)^T.
\]

\begin{lemma}\label{lem:direct-sum}
Suppose that each \(A_j\) has a detector \(R_j=K_j+t_jt_j^T\),
where \(K_j\in\{0,1,-1\}^{E_j\times E_j}\) is skew-symmetric and
\(K_j\bmod2=A_j+t_jt_j^T\).
Then \(A\) has a detector \(R=K+tt^T\), where
\(K\in\{0,1,-1\}^{E\times E}\) is skew-symmetric and
\(K\bmod2=A+tt^T\).
\end{lemma}

\begin{proof}
Define \(R\) in the block order \(E_1,\ldots,E_s\) by
\[
 R_{ij}=
 \begin{cases}
  R_i,&i=j,\\
  2t_it_j^T,&i<j,\\
  0,&i>j.
 \end{cases}
\]
Set \(K=R-tt^T\). Its blocks are
\[
 K_{ij}=
 \begin{cases}
  K_i,&i=j,\\
  t_it_j^T,&i<j,\\
  -t_it_j^T,&i>j.
 \end{cases}
\]
Since each \(K_i\) is skew-symmetric and \(K_{ji}=-K_{ij}^T\)
for \(i\ne j\), the matrix \(K\) is skew-symmetric.
Its entries belong to \(\{0,1,-1\}\), because the entries of
each \(K_i\) belong to this set and those of each \(t_i\) are \(0\) or \(1\).

On the diagonal blocks, the hypothesis gives
\(K_{ii}\bmod2=A_i+t_it_i^T\).
For \(i\ne j\), reduction modulo \(2\) gives
\(K_{ij}\bmod2=t_it_j^T\), which is the \((i,j)\)-block of \(A+tt^T\),
since \(A\) is block diagonal. Hence \(K\bmod2=A+tt^T\).

For \(X\subseteq E\), put \(X_j=X\cap E_j\).
The matrix \(R[X]\) is block upper triangular with diagonal blocks \(R_j[X_j]\).
Therefore
\[
 \det R[X]
 =\prod_{j=1}^s\det R_j[X_j]
 =\prod_{j=1}^s\det\nolimits_2 A_j[X_j]
 =\det\nolimits_2 A[X].
\]
The second equality uses that each \(R_j\) detects \(A_j\).
The last equality follows from the block diagonal form of \(A[X]\).
The formula also covers empty \(X_j\), whose principal determinants are \(1\).
Thus \(R=K+tt^T\) detects \(A\).
\end{proof}


\section{Normal forms in principal-minor fibres}\label{sec:fibres}

Within an irreducible principal-minor fibre, we show that a prescribed
symmetric part \(ss^T\) can be recovered by diagonal similarity. We use
this fact to normalize detectors on connected support.

For a real matrix on \(\{1,\ldots,m\}\), a cut is a partition
\(\{1,\ldots,m\}=X\sqcup Y\) with \(|X|,|Y|\ge2\) such that
both cross blocks have rank one. For an irreducible matrix, these
are the cuts of \cite[Definition~2.5]{CGGR}, since neither cross
block can be zero.
For a fixed cut \(X\sqcup Y\), factor the two cross blocks as
\(pq^T\) and \(uv^T\), with \(p,v\) indexed by \(X\) and \(q,u\)
indexed by \(Y\). The cut-transpose along this cut is
\cite[Definition~2.10]{CGGR}
\[
 B=\begin{pmatrix}P&pq^T\\uv^T&Q\end{pmatrix}
 \quad\longmapsto\quad
 \tau(B)=\begin{pmatrix}P&pu^T\\qv^T&Q^T\end{pmatrix}.
\]
Each rank-one factorization is unique up to multiplying its two
factors by reciprocal nonzero scalars.

\begin{lemma}\label{lem:general-fibre}
Let \(m\ge1\), let \(s\in\mathbb R^m\), and let
\(K\in\mathbb R^{m\times m}\) satisfy \(K^T=-K\).
Suppose that \(R=K+ss^T\) and \(N\in\mathbb R^{m\times m}\)
are irreducible and satisfy
\[
 \det R[X]=\det N[X]\qquad(X\subseteq\{1,\ldots,m\}).
\]
Then there exist an invertible real diagonal matrix \(D\)
and a real skew-symmetric matrix \(J\) such that
\[
 D^{-1}ND=J+ss^T.
\]
\end{lemma}

\begin{proof}
By \cite[Theorem~1.1 and Definition~2.10]{CGGR}, a sequence of
cut-transposes takes \(R\) to a matrix diagonally similar to \(N\)
or \(N^T\).
We show that these operations preserve diagonal similarity to
a matrix \(B\) satisfying \(B+B^T=2ss^T\).

First suppose that \(R_0+R_0^T=2ss^T\), and let
\(\{1,\ldots,m\}=X\sqcup Y\) be a cut of \(R_0\). Use the blocks
\(P,Q\) and factors \(p,q,u,v\) from the cut-transpose formula above.
All four factors are nonzero because the cross blocks have rank one.
The equation \(R_0+R_0^T=2ss^T\) gives
\[
 P+P^T=2s[X]s[X]^T,\qquad
 Q+Q^T=2s[Y]s[Y]^T,
\]
and
\[
 pq^T+vu^T=2s[X]s[Y]^T.
\]
The last matrix has rank at most one.
If both pairs \(p,v\) and \(q,u\) were linearly independent, then
\[
 pq^T+vu^T=
 \begin{pmatrix}p&v\end{pmatrix}
 \begin{pmatrix}q^T\\u^T\end{pmatrix}
\]
would have rank two: the first factor is injective on \(\mathbb R^2\),
and the second maps onto \(\mathbb R^2\).
Thus either \(u=\lambda q\) or \(v=\lambda p\)
for some nonzero \(\lambda\in\mathbb R\).

If \(u=\lambda q\), set \(F=\operatorname{diag}(I_X,\lambda^{-1}I_Y)\).
Then
\[
 F^{-1}\tau(R_0)F
 =\begin{pmatrix}P&pq^T\\uv^T&Q^T\end{pmatrix}.
\]
Its cross blocks agree with those of \(R_0\), and transposing \(Q\)
leaves \(Q+Q^T\) unchanged. Hence
\[
 F^{-1}\tau(R_0)F+\bigl(F^{-1}\tau(R_0)F\bigr)^T=2ss^T.
\]

If \(v=\lambda p\), set \(F=\operatorname{diag}(I_X,\lambda I_Y)\).
Now
\[
 F^{-1}\tau(R_0)F
 =\begin{pmatrix}P&\lambda pu^T\\qp^T&Q^T\end{pmatrix}.
\]
The sum of its upper right block and the transpose of its lower left block is
\[
 \lambda pu^T+pq^T=vu^T+pq^T=2s[X]s[Y]^T.
\]
Together with the identities for \(P+P^T\) and \(Q+Q^T\),
this again gives the required symmetric part.
Neither case requires \(s[X]\) or \(s[Y]\) to be nonzero.

Next let \(A_0=G^{-1}R_0G\), where
\(G=\operatorname{diag}(G_X,G_Y)\) is invertible diagonal.
Diagonal similarity preserves the ranks of the cross blocks,
so a cut of \(A_0\) is also a cut of \(R_0\).
The cross blocks of \(A_0\) have factorizations
\[
 (G_X^{-1}p)(G_Yq)^T,\qquad
 (G_Y^{-1}u)(G_Xv)^T.
\]
Using these factors in the cut-transpose formula gives
\[
 \tau(A_0)=E_0^{-1}\tau(R_0)E_0,\qquad
 E_0=\operatorname{diag}(G_X,G_Y^{-1}).
\]
Thus the preceding calculation also applies when the matrix
before the cut-transpose is diagonally similar to \(R_0\).

This conclusion is independent of the chosen factorizations.
Indeed, replacing
\[
 (p,q,u,v)\quad\text{by}\quad
 (\alpha p,\alpha^{-1}q,\beta u,\beta^{-1}v)
\]
for nonzero \(\alpha,\beta\in\mathbb R\) gives the cut-transpose
\[
 \begin{pmatrix}P&\alpha\beta pu^T\\
    (\alpha\beta)^{-1}qv^T&Q^T\end{pmatrix}
 =H^{-1}\tau(R_0)H,\qquad
 H=\operatorname{diag}(I_X,\alpha\beta I_Y).
\]
Every factorization therefore yields a diagonally similar cut-transpose.

Diagonal similarity itself preserves the required property by composition.
Whole-matrix transpose preserves it as well: if
\[
 G^{-1}A_0G=J+ss^T,\qquad J^T=-J,
\]
then, since \(G\) is diagonal,
\[
 GA_0^TG^{-1}=-J+ss^T.
\]

Starting from \(R+R^T=2ss^T\), apply the cut-transpose calculation
successively along the sequence. Its last matrix has the required
property, so diagonal similarity gives the property for \(N\) or
\(N^T\). In the latter case, apply the transpose calculation to
obtain it for \(N\).
There is an invertible real diagonal matrix \(D\) such that
\(B=D^{-1}ND\) satisfies \(B+B^T=2ss^T\).
Set \(J=B-ss^T\). Then
\[
 J^T=B^T-ss^T=2ss^T-B-ss^T=-J,
\]
and \(D^{-1}ND=J+ss^T\).
\end{proof}

\begin{lemma}\label{lem:general-local-normal}
Let \(B=B^T\in\Ftwo^{X\times X}\), where \(X\) is finite and nonempty,
and put \(s=\operatorname{diag}B\).
Suppose that the off-diagonal support graph of \(B\) is connected
and that \(B\) has a detector \(R=K+ss^T\) with
\(K\in\mathbb R^{X\times X}\) skew-symmetric.
Then every real detector \(N\) of \(B\) is diagonally similar to
\(J+ss^T\) for some skew-symmetric \(J\in\{0,1,-1\}^{X\times X}\).
\end{lemma}

\begin{proof}
Let \(M\) be either \(N\) or \(R\).
For each support edge \(ij\), equation~\eqref{eq:detector-pairs}
gives \(M_{ij}M_{ji}=2s_is_j-1\in\{1,-1\}\), so both entries are nonzero.

Any path in the support graph therefore gives a directed path in either
direction in the directed graph of nonzero off-diagonal entries of \(M\).
Since the support graph is connected, the directed graph is strongly connected.
When \(|X|=1\), it is strongly connected by the one-vertex convention.
Thus both \(N\) and \(R\) are irreducible.
They also have the same principal minors, since
\[
 \det N[Y]=\det\nolimits_2 B[Y]=\det R[Y]\qquad(Y\subseteq X).
\]
Lemma~\ref{lem:general-fibre} therefore gives an invertible real diagonal
matrix \(D\) and a real skew-symmetric matrix \(J\) such that
\[
 C:=D^{-1}ND=J+ss^T.
\]
Diagonal similarity preserves principal minors, so \(C\) detects \(B\).
For distinct \(i,j\in X\), the two-index calculation
\eqref{eq:rankone-pair} gives
\[
 J_{ij}^2=\det C[\{i,j\}]=\det\nolimits_2 B[\{i,j\}]\in\{0,1\}.
\]
The entries of \(J\) are real, so \(J_{ij}\in\{0,1,-1\}\)
for \(i\ne j\). Its diagonal entries are zero by skew-symmetry.
Consequently \(J\in\{0,1,-1\}^{X\times X}\).
\end{proof}


\section{Normalization in the zero-diagonal case}\label{sec:zero-diagonal}

\begin{lemma}\label{lem:zero-diagonal}
Let \(E\) be finite, let \(A=A^T\in\Ftwo^{E\times E}\) have zero
diagonal, and suppose that \(M\in\mathbb R^{E\times E}\) detects \(A\).
Then \(A\) has a PU signing.
If \(E\ne\varnothing\) and the off-diagonal support graph \(G\) of
\(A\) is connected, there is a positive invertible diagonal matrix
\(D\) such that \(K=D^{-1}MD\) is such a signing.
\end{lemma}

\begin{proof}
The empty matrix supplies the assertion for \(E=\varnothing\).
We first treat connected \(G\). For one vertex the singleton
equation gives \(M=[0]\), so take \(D=[1]\).
For larger \(G\), equation~\eqref{eq:detector-pairs} reads
\begin{equation}\label{eq:zero-pairs}
 M_{ii}=0,\qquad M_{ij}M_{ji}=-A_{ij}\quad(i\ne j).
\end{equation}
In particular, both entries on a support edge are nonzero.

\emph{Entries corresponding to nonedges.}
We prove \(M_{ij}=M_{ji}=0\) for every nonedge \(ij\), by induction
on the distance \(d\ge2\) between its endpoints in \(G\).
Choose a shortest path \(v_0=i,v_1,\ldots,v_d=j\).
It is induced, since a chord would shorten it.
Every nonconsecutive pair other than \(v_0,v_d\) has graph
distance at most its separation on this path, and that separation
is smaller than \(d\). The induction hypothesis eliminates both
entries of each such pair. For \(d=2\) there are no pairs of this
kind, so the same argument starts the induction.

Put \(S=\{v_0,\ldots,v_d\}\), \(u=M_{v_0v_d}\), and
\(v=M_{v_dv_0}\). Equation~\eqref{eq:zero-pairs} gives \(uv=0\).
Let \(P\) be the matrix obtained from \(M[S]\) by setting these
two corner entries to zero. Its nonzero entries lie on consecutive
path pairs, and its diagonal is zero.
Expansion along the first row and then along the column containing
the reverse first-edge entry gives
\[
 \det P=-M_{v_0v_1}M_{v_1v_0}
            \det P[\{v_2,\ldots,v_d\}]
        =\det P[\{v_2,\ldots,v_d\}].
\]
Repeating this calculation, with determinant \(1\) on the empty
set and \(0\) on a singleton, gives \(\det P=1\) when \(|S|\) is
even and \(\det P=0\) when \(|S|\) is odd.
The binary adjacency matrix of the path has the same recurrence
and hence the same integer \(0/1\) determinant.
Since \(M\) detects \(A\), it follows that \(\det M[S]=\det P\).

In the permutation expansion of \(\det M[S]\), a nonzero term
using the corner \(v_0\to v_d\) but not its reverse must return
from \(v_d\) to \(v_0\) along the entire path.
It cannot return through fewer vertices because there are no
other nonconsecutive entries. The reverse corner gives the
oppositely directed cycle. Terms using both corners contain
the factor \(uv=0\). All remaining terms occur in \(\det P\).
Thus
\begin{equation}\label{eq:zero-path}
 0=(-1)^d\left(
 u\prod_{a=0}^{d-1}M_{v_{a+1}v_a}
 +v\prod_{a=0}^{d-1}M_{v_av_{a+1}}
 \right).
\end{equation}
Both products are nonzero by \eqref{eq:zero-pairs}.
If \(u\ne0\), then \(v=0\), contradicting
\eqref{eq:zero-path}; if \(v\ne0\), the same equation gives a
contradiction with \(u=0\). Hence \(u=v=0\).
The induction proves that both entries corresponding to each nonedge of \(G\) vanish.

\emph{Absolute products around cycles.}
Let \(C=(v_1,\ldots,v_k,v_1)\) be an induced cycle of \(G\), with \(k\ge3\), and put
\[
 P_C=\prod_{a=1}^k M_{v_av_{a+1}},\qquad v_{k+1}=v_1.
\]
Multiplying \eqref{eq:zero-pairs} on its edges shows that the
reverse product is \((-1)^k/P_C\).
In the cycle restriction, the preceding nonedge argument has
eliminated every entry outside the cycle edges.
For odd \(k\), the only nonzero determinant terms are the two
directed full cycles. For even \(k\), there are also the two
perfect matchings. Each gives a permutation consisting of transpositions,
with factor \(-M_{ij}M_{ji}=1\) for each matched edge.
The full-cycle permutation has sign \((-1)^{k-1}\). Therefore
\begin{equation}\label{eq:zero-cycle}
 \det M[V(C)]=
 \begin{cases}
  P_C-P_C^{-1},&k\text{ odd},\\
  2-P_C-P_C^{-1},&k\text{ even}.
 \end{cases}
\end{equation}
The binary cycle adjacency matrix annihilates the nonzero all-one
column, since each row has two ones. Its determinant is zero.
Detection sets \eqref{eq:zero-cycle} equal to zero.
For odd \(k\) this gives \(P_C^2=1\).
For even \(k\), multiplying by \(P_C\ne0\) gives
\((P_C-1)^2=0\). In both cases,
\begin{equation}\label{eq:zero-absolute-balance}
 \prod_{a=1}^k |M_{v_av_{a+1}}|=1.
\end{equation}

Equation~\eqref{eq:zero-absolute-balance} holds for every simple
cycle. To prove this, induct on cycle length.
An induced cycle was just treated. A chord divides any other
simple cycle into two shorter simple cycles.
Traverse them in the original cycle direction and in opposite
directions on the chord. Multiplying their equations cancels
the chord because \(|M_{ij}||M_{ji}|=1\).
The remaining product is the product for the original cycle.

\emph{Diagonal normalization.}
Choose a rooted spanning tree of \(G\), set \(d_{\rm root}=1\),
and, along each parent-to-child edge \(i\to j\), set
\[
 d_j=d_i/|M_{ij}|.
\]
All \(d_i\) are positive.
For \(D=\operatorname{diag}(d_i:i\in E)\), each tree-edge entry
of \(D^{-1}MD\) has absolute value \(1\).
Along any simple cycle, the factors \(d_j/d_i\) in the entries of
\(D^{-1}MD\) cancel, so the absolute product remains \(1\) by
\eqref{eq:zero-absolute-balance}. In the cycle formed by a non-tree
edge and its tree path, all tree-edge factors are \(1\), including
those traversed in reverse by \eqref{eq:zero-pairs}. The non-tree
factor must also be \(1\).
The product of the two transformed entries on an edge is still
\(-1\), by \eqref{eq:zero-pairs}, so those entries are opposite
signs. The diagonal entries and the entries corresponding to nonedges of \(G\) are zero.
Thus \(K=D^{-1}MD\) is skew-symmetric, has entries in
\(\{0,1,-1\}\), and reduces to \(A\) modulo \(2\).
Its principal minors equal those of \(M\), so it is PU.

Finally, for disconnected \(G\), apply the connected argument to
the restriction \(M[C]\) on each component \(C\), and take the
block direct sum of the resulting signings.
For \(X\subseteq E\), its determinant is the product of the
component determinants. The same product describes the binary
determinant of \(A[X]\), because \(A\) is block diagonal.
Each factor is \(0\) or \(1\), so this is an integer detection
identity. The direct sum is a PU signing of \(A\).
\end{proof}


\section{Proof of the main theorem}\label{sec:general-binary}
\label{sec:general-bouquet}

The implication from real detectability to PU-orientability uses a
minimal counterexample. Proper restrictions constrain its support
graph, and feasible singleton pivots exclude the remaining graphs.

\subsection{Non-cut vertices and small matrices}
For a graph \(G\), write \(V(G)\) for its vertex set, \(G[X]\) for
the subgraph induced by \(X\subseteq V(G)\), and \(G-v\) for
\(G[V(G)\setminus\{v\}]\).
A cut vertex of a connected graph is a vertex whose deletion
disconnects the graph. Its blocks are its maximal 2-connected
subgraphs, with each bridge also counted as a two-vertex block.
The block-cut tree has a node for each block and each cut vertex;
a block node is adjacent to the cut vertices it contains.
For a connected graph with at least two vertices this incidence
graph is a tree \cite[Lemma~3.1.4]{Diestel}.

\begin{lemma}\label{lem:general-block-count}
Let \(G\) be a finite connected simple graph with at least two vertices.
Let \(k\) be its number of non-cut vertices, and let \(\mathcal B\)
be its block-cut tree, counting each bridge as a two-vertex block.
Then
\begin{equation}\label{eq:block-count}
 k=2+\sum_{\text{blocks }B}(|V(B)|-2)
       +\sum_{\text{cut vertices }v}(\deg_{\mathcal B}(v)-2).
\end{equation}
If \(G\) is not a path, then \(k\ge3\).
If \(k=3\), then \(G\) is a subdivision of \(K_{1,3}\), or a triangle
with a path attached at each vertex. In the latter case the paths
are otherwise disjoint, and each may have length zero.
\end{lemma}

\begin{proof}
Write \(b,c\) for the numbers of blocks and cut vertices and put
\(S=\sum_v\deg_{\mathcal B}(v)\), over cut vertices.
The tree has \(b+c-1\) edges, each incident with exactly one cut-vertex
node, so \(S=b+c-1\).
Each non-cut vertex belongs to one block; a cut vertex belongs to
as many blocks as its tree degree. Hence
\(\sum_B|V(B)|=k+S\).
The right side of~\eqref{eq:block-count} is
\(2+(k+S-2b)+(S-2c)=k\).
Every block has at least two vertices, and every cut-vertex node
has tree degree at least two. Thus all summands are nonnegative.
For \(k=2\), every block is a bridge
and every cut-vertex tree degree is two, giving a path.
For \(k=3\), the excess sum is one.
If all blocks are bridges, \(G\) is a tree, and each cut vertex has
degree in \(G\) equal to its degree in \(\mathcal B\). Exactly one
vertex has degree three and all others have degree at most two.
The three branches at the degree-three vertex are paths, so \(G\)
is a subdivision of \(K_{1,3}\).
Otherwise, there is one three-vertex block, necessarily a triangle,
and all other blocks are bridges. Each cut vertex belongs to two
blocks, so each triangle vertex has at most one attached branch and
the bridges in that branch form a path. The branches are disjoint
outside the triangle because \(\mathcal B\) is a tree.
\end{proof}

\begin{lemma}\label{lem:general-small}
If \(A\) is any binary symmetric matrix of order at most three,
then \(A^\sharp\) admits a PU-orientation.
\end{lemma}

\begin{proof}
The alternating binary matrix \(A^\sharp\) has order at most four.
Let \(L\) be obtained by replacing each nonzero upper-triangular
entry of \(A^\sharp\) by \(1\) or \(-1\), and then extending
skew-symmetrically. Thus \(L\bmod2=A^\sharp\). Principal minors
of orders zero and two are
already \(1\) or \(0\), and odd-order minors vanish.
Only order four needs consideration.
On four indices, its Pfaffian is
\[
 L_{12}L_{34}-L_{13}L_{24}+L_{14}L_{23}.
\]
Let \(r\in\{0,1,2,3\}\) count the nonzero terms.
Different terms use disjoint pairs of edges, so their signs can
be chosen independently. For \(r=0,1,2,3\), choose their sum to be,
respectively, \(0,1,0,1\).
The full determinant, the square of this Pfaffian, belongs to
\(\{0,1\}\). This gives the required PU-orientation.
\end{proof}

\subsection{Proof of the equivalence}

\begin{proof}[Proof of Theorem~\ref{thm:binary-detector}]
For \(E=\varnothing\), the empty matrices \(M,K\) satisfy
conditions 1 and 2, and
\(\widetilde{A^\sharp}=A^\sharp=[0]\) satisfies condition 3.
Assume \(E\ne\varnothing\).
For a detector as in condition 2, equation~\eqref{eq:special-reduction}
gives \(K\bmod2=A+tt^T\).
The implications \(2\Rightarrow1\) and \(2\Rightarrow3\) follow from
the definition and Lemma~\ref{lem:general-border}, respectively.
For \(3\Rightarrow2\), choose a PU-orientation
\(\widetilde{A^\sharp}\) of \(A^\sharp\). For each \(i\) with
\(t_i=1\), multiply row \(i\) and column \(i\) by
\((\widetilde{A^\sharp})_{ip}\in\{1,-1\}\).
Let \(S\) be the diagonal matrix of these signs, with \(S_{ii}=1\)
on the remaining indices, including \(p\). Conjugation by \(S=S^{-1}\)
preserves each principal minor by diagonal similarity. It also
preserves the reduction modulo \(2\) because \(S\equiv I\pmod2\).
The new border column restricted to \(E\) is \(t\). Write
\[
 \widetilde{A^\sharp}=
 \begin{pmatrix}K&t\\-t^T&0\end{pmatrix}.
\]
Equation~\eqref{eq:border-parity} gives
\[
 \det\bigl((K+tt^T)[X]\bigr)=\det\bigl(\widetilde{A^\sharp}[\alpha(X)]\bigr)
 \qquad(X\subseteq E).
\]
These minors belong to \(\{0,1\}\) because the border matrix is a
PU-orientation. Its upper left block satisfies
\(K\bmod2=A+tt^T\), so \((K+tt^T)\bmod2=A\).
The determinant congruence and the \(0/1\) values give
\(\det((K+tt^T)[X])=\det\nolimits_2 A[X]\).
It remains to prove \(1\Rightarrow3\).

\emph{Connected support and proper restrictions.}
Choose a binary symmetric matrix \(A\) of minimum order \(n=|E|\)
that has a real detector \(M\) and whose lift \(A^\sharp\) has no
PU-orientation.
Lemma~\ref{lem:general-small} gives \(n\ge4\).
If \(t=0\), Lemma~\ref{lem:zero-diagonal} supplies a PU signing of \(A\). Adding a zero row and column
signs \(A^\sharp=\operatorname{diag}(A,0)\), contradicting the
choice of \(A\). Hence \(t\ne0\).
Let \(G\) be the off-diagonal support graph of \(A\).
It is connected: otherwise each component is a smaller detectable
binary matrix, so minimality and \(3\Rightarrow2\) provide special detectors
on all components. Lemmas~\ref{lem:direct-sum} and
\ref{lem:general-border} then give a PU lift of \(A\).

For every nonempty proper \(X\subsetneq E\) with \(G[X]\) connected,
the detector equations on subsets of \(X\) show that \(M[X]\)
detects \(A[X]\). Minimality and \(3\Rightarrow2\) give a special
detector of \(A[X]\).
Lemma~\ref{lem:general-local-normal} gives an invertible real diagonal
matrix \(D_X\) and a skew-symmetric matrix
\(K_X\in\{0,1,-1\}^{X\times X}\) such that
\begin{equation}\label{eq:local-normalization}
 D_X^{-1}M[X]D_X=K_X+t[X]t[X]^T,\qquad K_X^T=-K_X.
\end{equation}
The full \(M\) cannot be diagonally similar to \(K+tt^T\) with
\(K^T=-K\): order-two detection would give \(K_{ij}^2\in\{0,1\}\),
and \eqref{eq:special-reduction} would give \(K\bmod2=A+tt^T\).
Lemma~\ref{lem:general-border} would then give a PU-orientation of
\(A^\sharp\).
By Lemma~\ref{lem:pivot-invariance}, every feasible singleton pivot
gives another counterexample on \(E\), so the same minimality argument
applies to the pivoted matrix.

\emph{Normalizing absolute values on support edges.}
First suppose \(G\) is not a chordless cycle.
On each support edge \(ij\), equation~\eqref{eq:detector-pairs} gives
\(|M_{ij}||M_{ji}|=1\).
The entries of a special detector on every support edge have
absolute value one. If \(t_i=t_j=1\), the target pair minor is
zero, so its skew entry is zero and both detector entries are \(1\).
If \(t_it_j=0\), the target pair minor is \(1\), so the skew
entry has square \(1\) and the rank-one entry \(t_it_j\) is zero.
For a support edge \(ij\) of a proper connected restriction \(G[X]\),
equation~\eqref{eq:local-normalization} gives
\[
 |M_{ij}|=\frac{|(D_X)_{ii}|}{|(D_X)_{jj}|}.
\]
Multiplying around a simple cycle with proper vertex set cancels the
ratios, so the product of the entry magnitudes around it equals \(1\).
Every spanning simple cycle has a chord, since \(G\) is not a
chordless cycle. The chord divides it into two cycles with proper
vertex sets. Multiplying their cycle products cancels the chord by
\(|M_{ij}||M_{ji}|=1\). Thus every simple cycle has product \(1\).
Choose a rooted spanning tree, put \(d_{\rm root}=1\), and set
\(d_j=d_i/|M_{ij}|\) on each tree edge directed away from the root.
Let \(D=\operatorname{diag}(d_i:i\in E)\) and \(M'=D^{-1}MD\).
On a tree edge directed from \(i\) to \(j\), the recursion gives
\(|M'_{ij}|=1\), and \(|M'_{ij}||M'_{ji}|=|M_{ij}||M_{ji}|=1\)
gives \(|M'_{ji}|=1\). Along any simple cycle, the factors \(d_j/d_i\)
cancel, so the product of the magnitudes of the entries of \(M'\)
along the cycle is still \(1\). For a non-tree edge, use the cycle
formed by this edge and its tree path. Every factor except the one
on the non-tree edge equals \(1\), so that factor is also \(1\).
For each \(X\) in \eqref{eq:local-normalization}, the new local
normalizer is \(D[X]^{-1}D_X\), since
\[
 (D[X]^{-1}D_X)^{-1}M'[X](D[X]^{-1}D_X)
 = D_X^{-1}M[X]D_X.
\]
We now use \(M\) and \(D_X\) to denote the normalized matrix and the
adjusted local normalizers, respectively.

\emph{Sign products around triangles.}
Suppose further that \(G\) is not a path.
Every pair of vertices lies in a proper connected induced subgraph:
take a shortest path between them; it is induced, and if it spans
all vertices the whole graph is a path.
On a support edge \(ij\) of \(G[X]\), the entries of both the
current \(M[X]\) and its local normal form have absolute value \(1\).
Thus \(|(D_X)_{jj}/(D_X)_{ii}|=1\) by \eqref{eq:local-normalization}.
Since \(G[X]\) is connected, all diagonal entries of \(D_X\) have the
same absolute value.
Write \(D_X\) as a scalar times
\(S_X=\operatorname{diag}(\sigma_i^{(X)}:i\in X)\), with
\(\sigma_i^{(X)}\in\{1,-1\}\). Adding \eqref{eq:local-normalization}
to its transpose gives
\[
 S_X(M[X]+M[X]^T)S_X=2t[X]t[X]^T,
 \qquad M_{ij}+M_{ji}=2\sigma_i^{(X)}\sigma_j^{(X)}t_it_j
 \quad(i,j\in X).
\]
Every pair belongs to such an \(X\), so \(M_{ij}+M_{ji}=0\) if
\(t_it_j=0\). Put \(T=\{i\in E:t_i=1\}\). For distinct \(i,j\in T\),
set \(q_{ij}=(M_{ij}+M_{ji})/2\in\{1,-1\}\).
If all triangles of the complete signed graph on \(T\) have product
one, choose \(a\in T\), which is nonempty since \(t\ne0\),
set \(\sigma_a=1,\ \sigma_i=q_{ai}\) on \(T\setminus\{a\}\),
and choose arbitrary signs off \(T\).
For distinct \(i,j\ne a\), the triangle equation on \(a,i,j\) gives
\(q_{ij}=q_{ai}q_{aj}=\sigma_i\sigma_j\); pairs containing \(a\)
satisfy the same identity by construction.
Set \(S=\operatorname{diag}(\sigma_i:i\in E)\) and \(B=S^{-1}MS\).
The off-diagonal identities give \(B_{ij}+B_{ji}=2t_it_j\) for
\(i\ne j\). On the diagonal, \(B_{ii}=M_{ii}=t_i=t_i^2\).
Hence \(B+B^T=2tt^T\), so \(K=B-tt^T\) is skew-symmetric.
This contradicts the exclusion of such a diagonal similarity above.
There must therefore be distinct \(a,b,c\in T\) with
\begin{equation}\label{eq:negative-triangle}
 q_{ab}q_{bc}q_{ca}=-1.
\end{equation}
No proper connected induced subgraph contains these three vertices,
since the three edge signs from that set would multiply to \(1\),
each local vertex sign appearing twice.
Every non-cut vertex consequently lies in \(\{a,b,c\}\):
deleting any other non-cut vertex would leave such a subgraph.
Lemma~\ref{lem:general-block-count} shows there are exactly three
non-cut vertices, namely \(a,b,c\), and classifies \(G\) as a
subdivision of \(K_{1,3}\) or a triangle with three attached paths.

There is no fourth vertex \(v\) with \(t_v=1\).
The graphs \(G-a,G-b,G-c\) are connected, so
\eqref{eq:local-normalization} gives
triangle product one on \(vbc,vac,vab\), respectively. Multiplying gives
\[
 1=(q_{va}q_{vb}q_{vc})^2q_{ab}q_{bc}q_{ca}=q_{ab}q_{bc}q_{ca},
\]
contrary to \eqref{eq:negative-triangle}.
Thus \(T=\{a,b,c\}\).

\emph{Pivoting at a leaf.}
The classified graph has a leaf \(\ell\).
In the subdivided claw this is an end of one of its three arms.
In the triangle case, \(n\ge4\) forces at least one attached
path to have positive length, and its far end is a leaf.
Its neighbour \(u\) is a cut vertex: deleting \(u\) isolates
\(\ell\) from the other vertices.
The leaf \(\ell\) is a non-cut vertex, so \(\ell\in T=\{a,b,c\}\).
Its neighbour \(u\) is a cut vertex, so \(u\notin T\).
Thus \(t_\ell=1\) and \(t_u=0\).
Put \(A'=A*\{\ell\}\). The binary singleton formula gives
\[
 A'_{ij}=A_{ij}+A_{i\ell}A_{\ell j}\quad(i,j\ne\ell),\qquad
 A'_{\ell j}=A_{\ell j}\quad(j\ne\ell),\qquad A'_{\ell\ell}=1.
\]
Because \(\ell\) has only neighbour \(u\), the only changed entry is
\(A'_{uu}=A_{uu}+A_{u\ell}A_{\ell u}=0+1=1\). In particular, the
support graph is still \(G\).
The new minimal counterexample has four diagonal entries equal to \(1\) on the
same support graph, contradicting the conclusion just proved for every
minimal counterexample whose support is neither a path nor a chordless cycle.
Thus every minimal counterexample has support a path or a
chordless cycle.

\emph{Eliminating paths and chordless cycles.}
If \(G\) is a chordless cycle, choose \(e\) with \(t_e=1\).
Since \(n\ge4\), its two neighbours are nonadjacent.
Pivoting at \(e\) adds the edge between them and preserves every
original edge. The new graph is connected and has two degree-three
vertices, so is neither a path nor a chordless cycle.
If \(G\) is a path with an internal vertex \(e\) satisfying
\(t_e=1\), pivot there.
This adds the edge between its two neighbours, retaining the
Since \(n\ge4\), at least one of the two neighbours also has a
neighbour farther along the original path. Its degree increases from
two to three, so the new graph is neither a path nor a chordless cycle.
Finally, if every internal vertex has diagonal entry zero, \(t\ne0\)
gives an endpoint with diagonal entry \(1\). Pivoting there preserves the
path and changes its internal neighbour's diagonal entry from \(0\) to
\(1\). A pivot at that neighbour is then feasible and gives the preceding case.
All cases contradict the restriction on every minimal counterexample.
No counterexample exists.
\end{proof}

\begin{corollary}\label{cor:binary-regularity}
Let \(E\) be finite and \(A=A^T\in\Ftwo^{E\times E}\).
Then \(A\) has a real detector if and only if
\(\widehat{D(A)}\) is regular.
\end{corollary}

\begin{proof}
By \eqref{eq:lift-minors}, \(D(A^\sharp)=\widehat{D(A)}\).
Theorem~\ref{thm:binary-detector} identifies real detectability
of \(A\) with PU-orientability of \(A^\sharp\), and
Lemma~\ref{lem:regular-signing} identifies the latter condition
with regularity of \(D(A^\sharp)\).
\end{proof}

\begin{corollary}\label{cor:general-connected-normal}
Let \(E\) be finite and nonempty, let \(A=A^T\in\Ftwo^{E\times E}\),
and put \(t=\operatorname{diag}A\). Suppose that the off-diagonal
support graph of \(A\) is connected. Every real detector \(M\)
of \(A\) is diagonally similar to \(K+tt^T\), where \(K\) is integer
skew-symmetric with entries in \(\{0,1,-1\}\).
\end{corollary}

\begin{proof}
Fix a real detector \(M\) of \(A\). Theorem~\ref{thm:binary-detector}
provides a special detector of \(A\), so
Lemma~\ref{lem:general-local-normal} applies to \(M\) on \(E\).
\end{proof}

\begin{example}\label{ex:disconnected-normal}
Let
\[
 A=\begin{pmatrix}0&0\\0&0\end{pmatrix}\quad\text{over }\Ftwo,
 \qquad M=\begin{pmatrix}0&1\\0&0\end{pmatrix}\quad\text{over }\mathbb R.
\]
Every nonempty principal minor of each matrix is zero, so \(M\)
detects \(A\). Here \(t=0\).
A special detector must be skew-symmetric, and its pair determinant
is \(K_{12}^2=0\); hence it is the zero matrix.
For an invertible diagonal \(D=\operatorname{diag}(d_1,d_2)\),
the \((1,2)\) entry of \(D^{-1}MD\) is \(d_2/d_1\ne0\).
Thus this given detector cannot be diagonally similar to a special
detector. The existence conclusion of
Theorem~\ref{thm:binary-detector} still holds, with \(K=0\).
\end{example}


\section{Generating polynomials}\label{sec:generating-polynomials}

For \(A=A^T\in\Ftwo^{E\times E}\), let \(\mathbf x=(x_i)_{i\in E}\)
be a family of indeterminates and define \(Q_A\in\mathbb Z[x_i:i\in E]\) by
\[
 Q_A(\mathbf x)=\sum_{X\subseteq E}\det\nolimits_2 A[X]
 \prod_{i\in X}x_i.
\]
A polynomial is Hurwitz stable if it is nonzero at every complex
tuple whose coordinates have positive real part
\cite[Definition~4.1]{DK}.

\begin{corollary}\label{cor:hurwitz}
Let \(E\) be finite and \(A=A^T\in\Ftwo^{E\times E}\).
If \(A\) has a real detector, then \(Q_A\) is Hurwitz stable.
\end{corollary}

\begin{proof}
If \(E=\varnothing\), then \(Q_A=1\). Otherwise, Theorem~\ref{thm:binary-detector}
gives a special detector \(R=K+tt^T\). Since \(R\) detects \(A\),
\[
 Q_A(\mathbf x)=\sum_{X\subseteq E}\det R[X]\prod_{i\in X}x_i.
\]
Here \(K\) is real skew-symmetric and \(t\) is a real column vector, so
\cite[Lemma~4.4]{DK} gives Hurwitz stability.
\end{proof}

For a regular even delta-matroid, the proof of
\cite[Theorem~4.7]{DK} obtains negative real zeros after grouping
the nonzero size parity, and \cite[Corollary~4.9]{DK} gives the
resulting ultra-log-concavity. The next corollary applies that
principle to the parity lift, with an explicit factorization.

\begin{corollary}
\label{cor:paired-real-roots}
Let \(E\) be finite, let \(A=A^T\in\Ftwo^{E\times E}\), and
suppose that \(A\) has a real detector.
For each integer \(j\), define
\[
 q_j=\sum_{\substack{X\subseteq E\\|X|=j}}\det\nolimits_2 A[X],
\]
so \(q_j=0\) when \(j<0\) or \(j>|E|\).
For each integer \(i\ge0\), put \(c_i=q_{2i-1}+q_{2i}\), and define
\[
 C_A(z)=\sum_{i\ge0}c_i z^i.
\]
If \(r=\deg C_A\), there are positive real numbers
\(\lambda_1,\ldots,\lambda_r\) such that
\begin{equation}\label{eq:paired-factorization}
 C_A(z)=\prod_{j=1}^r(1+\lambda_j^2z).
\end{equation}
Thus every zero of \(C_A\) is real and strictly negative,
\(c_i>0\) for \(0\le i\le r\), and
\begin{equation}\label{eq:paired-newton}
 \left(\frac{c_i}{\binom ri}\right)^2
 \ge
 \frac{c_{i-1}}{\binom r{i-1}}\,
 \frac{c_{i+1}}{\binom r{i+1}}
 \qquad(1\le i<r).
\end{equation}
For \(r=0\), the product is \(1\) and the inequalities have
no indices.
\end{corollary}

\begin{proof}
If \(E=\varnothing\), then \(q_0=c_0=1\), all other coefficients
vanish, and \(C_A=1\). Assume \(E\ne\varnothing\).
By Theorem~\ref{thm:binary-detector}, choose a PU-orientation \(L\)
of \(A^\sharp\). Each principal minor of \(L\) belongs to \(\{0,1\}\).
Since \(L\bmod2=A^\sharp\), equation~\eqref{eq:lift-minors} gives
\[
 \det L[\alpha(X)]=\det\nolimits_2 A^\sharp[\alpha(X)]
 =\det\nolimits_2 A[X]\qquad(X\subseteq E),
\]
where the binary determinants are read as integers.
An even subset of \(E\cup\{p\}\) of size \(2i\) either avoids
\(p\), and is \(\alpha(X)\) for \(|X|=2i\), or contains \(p\),
and is \(\alpha(X)\) for \(|X|=2i-1\).
All odd principal minors of \(L\) are zero. Determinant
expansion therefore gives
\begin{equation}\label{eq:paired-characteristic}
 \det(I+sL)
 =\sum_{i\ge0}c_i s^{2i}=C_A(s^2).
\end{equation}

Let \(r_0=\operatorname{rank}(L)/2\).
By the real skew-symmetric canonical form
\cite[Sections~I.B--I.C]{Wimmer}, an orthogonal change of basis
puts \(L\) into a block diagonal matrix consisting of a zero block
and \(r_0\) blocks
\(\left(\begin{smallmatrix}0&\lambda_j\\-\lambda_j&0\end{smallmatrix}\right)\),
where \(\lambda_j>0\).
Taking the determinant block by block gives
\[
 \det(I+sL)=\prod_{j=1}^{r_0}(1+\lambda_j^2s^2).
\]
Comparing with \eqref{eq:paired-characteristic} and using injectivity
of \(z\mapsto s^2\) gives
\(C_A(z)=\prod_{j=1}^{r_0}(1+\lambda_j^2z)\).
The coefficient of \(z^{r_0}\) is \(\prod_{j=1}^{r_0}\lambda_j^2>0\),
so \(r=r_0\).
The factors show both the asserted zero locations and positivity
of every coefficient from degree \(0\) to degree \(r\).
Newton's inequalities \cite[Theorem~1.1]{Niculescu} for elementary
symmetric functions of the positive numbers
\(\lambda_1^2,\ldots,\lambda_r^2\) give
\eqref{eq:paired-newton}.
\end{proof}

Inequality~\eqref{eq:paired-newton} says that \((c_i)_{i=0}^r\)
is ultra-log-concave: after division by \(\binom ri\), the
sequence is log-concave. Its positivity gives no internal zeros.
Both corollaries also apply when \(\widehat{D(A)}\) is regular,
by Corollary~\ref{cor:binary-regularity}.
Ding and Kim obtain the corresponding ribbon-graph conclusions
on their pseudo-orientable class
\cite[Theorems~1.3--1.4 and Corollary~4.10]{DK}.


\section{Application to ribbon bouquets}\label{sec:bouquet-applications}

\begin{definition}
A ribbon graph is a surface with boundary formed from vertex
discs and edge bands, each band attached along two disjoint
boundary intervals. A bouquet is a ribbon graph with one vertex
disc. A loop is twisted if its band together with the vertex
disc is a M\"obius band, and untwisted if it is an annulus.
A quasi-tree is a spanning ribbon subgraph with one boundary
component. These ribbon-graph conventions follow
\cite[Section~2.1]{DK}.
\end{definition}

The attachment intervals of a bouquet give a chord diagram on
the boundary circle of its vertex. Two distinct loops interlace
when their endpoints alternate on that circle.
The framed intersection matrix \(A\) is symmetric over \(\Ftwo\),
with \(A_{ii}=1\) for twisted loops and \(A_{ij}=1\) for
interlacing distinct loops.
A graph obtained as the interlacement graph of chords is a
circle graph; see \cite[Definition~2.4 and Proposition~2.21]{DK}
for the chord-diagram and matrix description.

For a bouquet, the quasi-tree edge sets are precisely the
subsets \(X\) for which \(A[X]\) is nonsingular
\cite[Proposition~2.21]{DK}. Thus \(D(A)\) is its ribbon-graph
delta-matroid. We say that a real matrix detects the bouquet
when its principal minor on \(X\) is \(1\) for a quasi-tree
edge set and \(0\) for every other edge set.

\begin{corollary}[Whistler's conjecture]\label{cor:bouquet-main}
Let \(\mathbb B\) be a finite ribbon bouquet with edge set \(E\)
and framed intersection matrix \(A\).
Put \(t=\operatorname{diag}A\), and choose \(p\notin E\).
The following conditions are equivalent:
\begin{enumerate}
\item Some real \(E\times E\) matrix detects \(\mathbb B\).
\item Some detector of \(\mathbb B\) is \(K+tt^T\), where
 \(K\) is integer skew-symmetric with entries in \(\{0,1,-1\}\).
\item The binary matrix
 \(A^\sharp=\left(\begin{smallmatrix}A+tt^T&t\\t^T&0\end{smallmatrix}\right)\)
 admits a PU-orientation.
\item The parity lift of \(D(A)\) is regular.
\end{enumerate}
\end{corollary}

\begin{proof}
The principal-minor description of quasi-trees identifies
condition 1 with real detectability of \(A\), and condition 2
with existence of a special detector of \(A\).
Theorem~\ref{thm:binary-detector} gives the equivalence of
conditions 1--3.
By \eqref{eq:lift-minors}, \(D(A^\sharp)=\widehat{D(A)}\).
Lemma~\ref{lem:regular-signing} gives the equivalence of conditions 3 and 4.
\end{proof}

The implication from condition 1 to condition 3 is
\cite[Conjecture~6.2]{Whistler}.
For any bouquet satisfying these equivalent conditions,
Corollary~\ref{cor:hurwitz} applies to its multivariate
quasi-tree generating polynomial.
In Corollary~\ref{cor:paired-real-roots}, \(q_j\) is then the
number of its quasi-trees with \(j\) edges, so the polynomial
counting sizes \(2i-1\) and \(2i\) has only negative real zeros.


\section{A converse question for Hurwitz stability}\label{sec:open-question}

\begin{question}
\label{ques:stability-converse}
Let \(E\) be finite and let \(A=A^T\in\Ftwo^{E\times E}\).
If \(Q_A\) is Hurwitz stable, does there exist a real matrix
\(M\in\mathbb R^{E\times E}\) such that
\[
 \det M[X]=\det\nolimits_2 A[X]
 \qquad\text{for every }X\subseteq E,
\]
\end{question}

By Theorem~\ref{thm:binary-detector} and Corollary~\ref{cor:binary-regularity},
the question is equivalent to asking whether Hurwitz stability of \(Q_A\)
implies PU-orientability of \(A^\sharp\), or equivalently regularity of
\(\widehat{D(A)}\).

When \(A\) has zero diagonal,
\(A^\sharp=\operatorname{diag}(A,0)\).
In this case, the question is whether Hurwitz stability
of \(Q_A\) implies that \(A\) has a PU signing.

\section*{Acknowledgements}
GPT-6.1 Sol and DeepSeek V4.1 Flash were used extensively to assist
with the development and preparation of this work.

\end{document}